\documentclass[12pt,a4paper]{amsart}

\usepackage{amsfonts, amsmath, amssymb, amsthm}
\theoremstyle{plain}
\newtheorem{thm}{Theorem}[section]
\newtheorem*{thm*}{Theorem}
\newtheorem{lem}[thm]{Lemma}
\newtheorem{prop}[thm]{Proposition}
\newtheorem{cor}[thm]{Corollary}

\theoremstyle{definition}
\newtheorem{defn}[thm]{Definition}
\newtheorem{exmp}[thm]{Example}

\theoremstyle{remark}

\usepackage{txfonts}
\usepackage{mathbbol}
\usepackage{ae}

\usepackage{color}
\usepackage{graphicx}
\usepackage{url}
\usepackage{overpic}
\usepackage{subfigure}
\usepackage{booktabs}
\usepackage{paralist}
\usepackage[a4paper]{anysize}

\newcommand\cF{{\mathcal F}}

\newcommand\NN{{\mathbb N}}
\newcommand\RR{{\mathbb R}}

\newcommand\QQ{{\mathbb Q}}

\newcommand\SetOf[2]{\left\{#1\vphantom{#2}\,\right.\left|\,\vphantom{#1}#2\right\}}
\newcommand\smallSetOf[2]{\{#1\,|\,#2\}}

\newcommand\conv{\operatorname{conv}}

\newcommand\bounded[1]{\partial_{\rm fin}{#1}}
\newcommand\interior[1]{\operatorname{int}{#1}}

\newcommand\outdeg{\operatorname{outdeg}}

\newcommand\argmax{\operatorname{argmax}}
\newcommand\dmax{d_{\max}}
\newcommand\dmin{d_{\min}}
\newcommand\Pmax{P_{\max}}
\newcommand\Tmax{T_{\max}}
\newcommand\Tmin{T_{\min}}
\newcommand\supp{\operatorname{supp}}
\newcommand\opt[1]{{{#1}_{\operatorname{opt}}}}

\renewcommand{\phi}{\varphi}

\newcommand\whyrelation[2]{\stackrel{\text{#1}}{#2}}

\graphicspath{{fig/} {jvx/}}

\begin{document}

\title{Bounds on the $f$-Vectors of Tight Spans}
\author[Herrmann \& Joswig]{Sven Herrmann and Michael Joswig}
\address{Sven Herrmann \& Michael Joswig, Fachbereich Mathematik, AG~7, TU Darmstadt, 64289 Darmstadt, Germany}
\email{sherrmann@mathematik.tu-darmstadt.de}
\email{joswig@mathematik.tu-darmstadt.de}
\thanks{The second author is partially supported by Deutsche Forschungsgemeinschaft, DFG Research Group
  ``Polyhedral Surfaces.''}
\date{\today}

\begin{abstract}
  The tight span $T_d$ of a metric $d$ on a finite set is the subcomplex of bounded faces of an unbounded polyhedron
  defined by~$d$.  If $d$ is generic then $T_d$ is known to be dual to a regular triangulation of a second hypersimplex.
  A tight upper and a partial lower bound for the face numbers of $T_d$ (or the dual regular triangulation) are
  presented.
\end{abstract}

\maketitle

\section{Introduction}
\noindent
Associated with a finite metric~$d:\{1,\dots,n\}\times\{1,\dots,n\}\to\RR$ is the unbounded polyhedron
\[ P_d \ = \ \SetOf{x\in\RR^n}{x_i+x_j\ge d(i,j)\text{ for all $i,j$}} \quad . \] Note that the condition ``for all
$i,j$'' includes the diagonal case $i=j$, implying that $P_d$ is contained in the positive orthant and thus pointed.
Following Dress~\cite{MR753872} we call the polytopal subcomplex $T_d$ formed of the bounded faces of~$P_d$ the
\emph{tight span} of~$M$; see also Bandelt and Dress~\cite{MR858908}.  In Isbell's paper~\cite{MR0182949} the same
object arises as the \emph{injective envelope} of~$d$.  The metric $d$ is said to be \emph{generic} if the
polyhedron~$P_d$ is simple.

Up to a minor technicality, the tight span $T_d$ is dual to a regular subdivision of the \emph{second hypersimplex}
\[ \Delta_{n,2} \ = \ \conv\SetOf{e_i+e_j}{1\le i<j\le n} \quad ,\]
and the tight spans for generic metrics correspond to regular triangulations.

The tight spans of metric spaces with at most six points have been classified by Dress~\cite{MR753872} and Sturmfels and
Yu~\cite{MR2097310}; see also De Loera, Sturmfels, and Thomas~\cite{MR1357285} for further details.
Develin~\cite{Develin} obtained sharp upper and lower bounds for the dimension of a tight span of a metric on a given
number of points.  The present paper can be seen as a refined analysis of Develin's paper.  Our main result is the
following.

\begin{thm*}
  The number of $k$-faces in a tight span of a metric on $n$ points  is at most
  \[ 2^{n-2k-1}\frac{n}{n-k}\binom{n-k}{k} \quad , \]
  and for each $n$ there is a metric $\dmax^n$ uniformly attaining this upper bound.
\end{thm*}

In particular, this result says that there are no $k$-faces for $k>\lfloor n/2\rfloor$, which is Develin's upper bound
on the dimension of a tight span.  Since the vertices of the tight span correspond to the facets of a hypersimplex
triangulation, and since further $\Delta_{n,2}$ admits an unimodular triangulation, this upper bound of $2^{n-1}$ for
the number of vertices of~$T_d$ is essentially the volume of $\Delta_{n,2}$.  In fact, the normalized volume of
$\Delta_{n,2}$ equals $2^{n-1}-n$, but this minor difference will be explained later.

The paper is organized as follows.  We start out with a section on the combinatorics of unbounded convex polyhedra.
Especially, we are concerned with the situation where such a polyhedron, say $P$, of dimension~$n$, is simple, that is,
each vertex is contained in exactly $n$ facets.  It then turns out that the $h$-vector of the simplicial ball which is
dual to the bounded subcomplex of~$P$ has an easy combinatorial interpretation using the vertex-edge graph of~$P$.  This
is based on ---and at the same time generalizes--- a result of Kalai~\cite{MR964396}.  Further, translating Develin's
result on the upper bound of the dimension of a tight span to the dual, says that a regular triangulation of a second
hypersimplex $\Delta_{n,2}$ does not have any interior faces of dimension up to $\lfloor (n-1)/2\rfloor-1$.  As a
variation of a concept studied by McMullen~\cite{MR2070631} and others we call such triangulations \emph{almost small
  face free} or \emph{asff}, for short.  The Dehn-Sommerville equations for the boundary then yield strong restrictions
for the $h$-vector of an asff simplicial ball.  Applying these techniques to the specific case of hypersimplex
triangulations leads to the desired result.  The final two sections focus on the construction of extremal metrics.  Here
the metric $\dmax^n$ is shown to uniformly attain the upper bound on the $f$-vector.  The situation turns out to be more
complicated as far as lower bounds are concerned.  The paper concludes with a lower bound for the number of faces of
maximal dimension of a tight span of dimension $\lceil n/3\rceil$, which is Develin's lower bound.  Further we construct
a metric $\dmin^n$ which attains this lower bound.  However, we do not have a tight lower bound for the number of faces
of smaller dimension.  Our analysis suggests that such a result might require to classify all possible $f$-vectors of
tight spans, a task beyond the scope of this paper.

\section{Combinatorics of Unbounded Polyhedra}
\noindent
A \emph{(convex) polyhedron} is the intersection of finitely many affine halfspaces in Euclidean space.  Equivalently, it
is the set of feasible solutions of a linear program.  A polyhedron $P$ is called \emph{pointed} if it does not contain
any affine line or, equivalently, its lineality space is trivial.  Further, $P$ is pointed if and only if it has at
least one vertex.  A \emph{(convex) polytope} is a bounded polyhedron.  For basic facts about polytopes and polyhedra
the reader may consult Ziegler~\cite{MR1311028}.

For a not necessarily pointed bounded polyhedron $P$ we denote the face poset by $\cF(P)$.  If $P$ is bounded then
$\cF(P)$ is a Eulerian lattice.  Two pointed polyhedra are called \emph{combinatorially equivalent} if their face posets
are isomorphic.

A polyhedron $P$ is pointed if and only if it is projectively equivalent to a polytope.  For this reason one can always
think of a pointed polyhedron $P$ as a polytope $P'$ with one face marked: the \emph{face at infinity}.  However, this
is not the only way to turn an unbounded polyhedron into a polytope: Take an affine halfspace $H^+$ which contains all
the vertices of $P$ and whose boundary hyperplane $H$ intersects all the unbounded edges.

\begin{lem}
  The combinatorial type of the polytope $\bar{P}=P\cap H^+$ only depends on the combinatorial type of $P$.
\end{lem}

\begin{proof}
  The vertices of $\bar{P}$ come in two kinds: Either they are vertices of~$P$ or they are intersections of rays of~$P$
  with the hyperplane~$H$.  The rays can be recognized in the face poset of the unbounded polyhedron~$P$ as those edges
  which contain only one vertex.  The claim now follows from the fact that the face lattice of the polytope $\bar{P}$ is
  atomic, that is, each face of~$\bar{P}$ is the join of vertices of~$\bar{P}$.
\end{proof}

We call $\bar{P}$ the \emph{closure} of~$P$.

The vertices and the bounded edges of a polyhedron~$P$ form an abstract graph which we denote by~$\Gamma(P)$.  Note that
in the unbounded case the rays (or unbounded edges) of~$P$ are not represented in~$\Gamma(P)$.

An $n$-dimensional pointed polyhedron $P$ is \emph{simple} if each vertex is contained in exactly $n$ facets.  Clearly,
simplicity is a combinatorial property.  If $P$ is bounded, that is, $P$ is a polytope, then it is simple if and only if
the graph $\Gamma(P)$ is $n$-regular.


\begin{prop}
  The pointed polyhedron $P$ is simple if and only if its closure $\bar{P}$ is.
\end{prop}

\begin{proof}
  If $P$ is a simple polyhedron, then $P$ is combinatorially equivalent to a polyhedron $Q$ which is the intersection of
  (facet defining) affine halfspaces in general position.  Without loss of generality we can choose an affine hyperplane
  $H$ which is in general position with respect to the facets of~$Q$ and which has the property that $H^+$ contains the
  vertices of~$P$.  Then $Q\cap H$ is simple, that is, $\Gamma(Q\cap H)$ is $(n-1)$-regular.  By construction each
  vertex of $Q\cap H$ is contained in exactly one unbounded edge of~$Q$.  This implies that the graph of the closure
  $\Gamma(Q\cap H^+)$ is $n$-regular, whence $\bar{Q}=Q\cap H^+$ is simple.  The reverse implication is trivial.
\end{proof}

\begin{prop}
  The combinatorial type of $\bar{P}$ is determined by the $2$-skeleton $\cF_{\le 2}(P)$.
\end{prop}

\begin{proof}
  The unbounded edges of~$P$ are exactly those edges which contain exactly one vertex each.  Hence $\cF_{\le 2}(P)$
  determines the vertices of the face $P\cap H$ in the closure $\bar{P}=P\cap H^+$.  The edges of $P\cap H$ correspond
  to the unbounded $2$-faces of~$P$, that is, those $2$-faces which contain two unbounded edges.  Altogether $\cF_{\le
    2}(P)$ determines the graph of the simple polytope $\bar{P}$.  A result of Blind and Mani~\cite{MR921106} then
  yields the claim.
\end{proof}


The \emph{bounded subcomplex} $\bounded{P}$ of an unbounded polyhedron $P$ is the polyhedral subcomplex of the boundary
$\partial P$ of~$P$ which is formed of the bounded faces.  Clearly, $\bounded{P}$ is contractible.  The graph
$\Gamma(P)$ is the $1$-skeleton of the bounded subcomplex.

Kalai's proof~\cite{MR964396} of the aforementioned result of Blind and Mani~\cite{MR921106} is based on a
characterization of the $h$-vector of a simple polytope in terms of acyclic orientations of its graph.  The remainder of
this section is devoted to explaining how this can be extended to bounded subcomplexes of unbounded polyhedra.

Consider an $n$-dimensional pointed polyhedron $P\subset\RR^n$ which is unbounded and a generic linear objective
function~$\alpha:\RR^n\to\RR$.  Let us assume that $\alpha$ is \emph{generic} on $\bar{P}=P\cap H^+$, that is, it is
$1$--$1$ on the vertices of~$\bar{P}$.  This way each edge of $P$, bounded or not, is a directed arc, say, with the
decrease of~$\alpha$.  Let us assume further that $\alpha$ is \emph{initial} with respect to $\bar{P}\cap H=P\cap H$,
that is, there are no arcs pointing towards the face $\bar{P}\cap H$ of $\bar{P}$.  In the language of linear
optimization, this means that the linear program $\max\smallSetOf{\alpha x}{x\in P}$ is unbounded and that the reverse
linear program
\[ \min\SetOf{\alpha x}{x\in P} \]
has a unique optimal vertex.

For each vertex $v$ of $\bar{P}$ let the \emph{out-degree} $\outdeg v$, with respect to~$\alpha$, be the number of edges
in $\bar{P}$ which are incident with~$v$ and directed away from~$v$.  For any subset $U$ of the vertices of~$\bar{P}$ we
let
\[ h_i(U)=\#\SetOf{v\in U}{\outdeg v=i} \; .\]

\begin{prop}\label{prop:f-from-h}
  We have
  \[
  f_k(\bounded{P})\ =\ \sum_{i=k}^n\binom{i}{k}h_i(P) \quad .
  \]
\end{prop}

\begin{proof}
  Each non-empty bounded face $F$ of $P$ has a unique $\alpha$-maximal vertex~$v=\argmax\alpha(F)$.  Conversely, $F$ is
  the unique face of~$P$ which is spanned by the edges in~$F$ which are incident with $v$.  This way
  $\binom{i}{k}h_i(P)$ counts those $k$-dimensional faces $F\le \bar{P}$ whose maximal vertex is not in $\bar{P}\cap H$ and
  which has $\outdeg\argmax\alpha(F)=i$.
\end{proof}

Later we will be interested in maximizing the $f$-vector of the bounded subcomplexes of certain unbounded polyhedra.
Because the binomial coefficients are non-negative, the previous proposition implies that maximizing the $f$-vector is
equivalent to maximizing the $h$-vector.

\section{Combinatorics of Simplicial Balls}\label{section:simplicial-balls}

\noindent
For an arbitrary $n$-dimensional simplicial complex $K$ with $f$-vector $f(K)$ we can define its \emph{$h$-vector} by
letting
\begin{equation}\label{eq:h-from-f}
  h_k(K)=\sum_{i=0}^k(-1)^{k-i}\binom{n+1-i}{n+1-k}f_{i-1}(K) \quad .
\end{equation}
Moreover, the \emph{$g$-vector} is set to $g_0(K)=1$ and $g_k(K)=h_k(K)-h_{k-1}(K)$ for $k\ge 1$.

As a consequence of the Euler equation, iteratively applied to intervals in the face lattice, we obtain the
\emph{Dehn-Sommerville relations}.

\begin{thm}\label{thm:DS}
  For each simplicial $(n-1)$-sphere~$S$ we have
  \[
  h_k(S)\ =\ h_{n-k}(S) \quad .
  \]
\end{thm}

As a further consequence the $f$-vectors (or $g$- or $h$-vectors) of a simplicial ball and its boundary are related.

\begin{thm}{\rm (McMullen and Walkup~\cite{MR0298557})}\label{thm:DS-ball} For each simplicial $(n-1)$-ball~$B$ we have
  \[
  g_k(\partial B)\ =\ h_k(B)-h_{n-k}(B) \quad .
  \]
\end{thm}

See also Billera and Bj\"orner~\cite{MR1730171} and McMullen~\cite[Corollary 2.6]{MR2070631}.

Let $\interior{B}$ be the set of interior faces of the ball~$B$.  Although $\interior{B}$ is not a polyhedral complex we
nonetheless write $f(\interior{B}):=f(B)-f(\partial B)$ for its $f$-vector.  Formally, we can also define the $h$-vector
of the interior faces of a ball by using the equation~\eqref{eq:h-from-f}.

\begin{prop}\label{prop:h-int}
  For each simplicial $(n-1)$-ball~$B$ we have
  \[  h_{n-k}(B) \ = \ h_k(\interior{B}) \quad . \]
\end{prop}

\begin{proof}
  \begin{align*}
    h_{n-k}(B)\ \whyrelation{\ref{thm:DS-ball}}{=}\
    & h_k(B)-g_k(\partial B)\ \whyrelation{\eqref{eq:h-from-f}}{=}\
    \sum_{i=0}^k(-1)^{k-i}\binom{n-i}{n-k}f_{i-1}(B) \\
    & - \left(\sum_{i=0}^k(-1)^{k-i}\binom{n-i-1}{n-k-1}f_{i-1}(\partial B)
    - \sum_{i=0}^{k-1}(-1)^{k-i-1}\binom{n-i-1}{n-k}f_{i-1}(\partial B)\right)\\
    =\ & \sum_{i=0}^k(-1)^{k-i}\binom{n-i}{n-k}\bigl(f_{i-1}(\interior{B})+f_{i-1}(\partial B)\bigr)
    \;-\; \sum_{i=0}^k(-1)^{k-i}\binom{n-i}{n-k}f_{i-1}(\partial B)\\
    =\ & \sum_{i=0}^k(-1)^{k-i}\binom{n-i}{n-k}f_{i-1}(\interior{B}) \ \whyrelation{\eqref{eq:h-from-f}}{=} \  h_k(\interior{B}) \quad .
  \end{align*}
\end{proof}

The following proposition is due to McMullen~\cite[Proposition 2.4c]{MR2070631}. We include its simple proof for the
sake of completeness.

\begin{prop}\label{prop:asff-h-vanishing}
  Let $B$ be a simplicial $(n-1)$-ball without any interior faces of dimension up to~$e$.  Then
  \[
  h_k(B)\ = \ 0 \; \text{ for $k\ge n-e-1$}\qquad  \text{and}\qquad h_k(B)\ =\ g_k(\partial B) \; \text{ for $k\le e+1$}\quad .
  \]
\end{prop}

\begin{proof}
  Our assumption on the interior faces says that $f_k(\interior{B})=0$ for $k\le e$.  From the proof of
  Proposition~\ref{prop:h-int} we see that
  \[
  h_{n-k}(B) \ = \ \sum_{i=0}^k(-1)^{k-i}\binom{n-i}{n-k}f_{i-1}(\interior{B}) \quad ,
  \]
  which directly proves $h_{n-k}(B)=0$ for $k\le e+1$.  Applying Theorem~\ref{thm:DS-ball} once again also proves the
  second claim.
\end{proof}

Of special interest is the case of a simplicial ball without small interior faces.  Following
McMullen~\cite[\S3]{MR2070631} we call a face $\sigma$ of a simplicial $(n-1)$-ball \emph{small} if
$\dim\sigma\le\lfloor(n-1)/2\rfloor$, and it is \emph{very small} if $\dim\sigma<\lfloor(n-1)/2\rfloor$.  A simplicial
$(n-1)$-ball is \emph{(almost) small-face-free}, abbreviated \emph{(a)sff}, if it does not have any (very) small
interior faces.

\begin{cor}\label{cor:defined-by-boundary-n-odd}
  The $f$-vector of an $(n-1)$-dimensional asff simplicial ball, for $n$ odd, is determined by the $f$-vector of its boundary.
\end{cor}

\begin{proof}
  Assume that $B$ is an $(n-1)$-dimensional asff simplicial ball.  Then we have
  \[
  f_k(B) \ = \ \sum_{i=k}^{n} \binom{i}{k}h_i(B)
  \ \whyrelation{\ref{prop:asff-h-vanishing}}{=} \
  \sum_{i=k}^{(n-1)/2} \binom{i}{k}g_i(\partial B) \quad .
  \]
\end{proof}

A similar computation shows the following analog for $n$ even.

\begin{cor}\label{cor:defined-by-boundary-n-even}
  The $f$-vector of an $(n-1)$-dimensional asff simplicial ball, for $n$ even, is determined by the $f$-vector of its
  boundary and $f_{n/2-1}=h_{n/2-1}$.
\end{cor}

A polytope is \emph{simplicial} if each proper face is a simplex.  Equivalently, its boundary complex is a simplicial
sphere.  In terms of cone polarity simplicity and simpliciality of polytopes are dual notions.  In this way, the bounded
subcomplex $\bounded{P}$ of an unbounded simple $n$-polyhedron $P$ becomes the set of interior faces of a simplicial
$(n-1)$-ball $B(P)$ in the boundary of the polar dual $\bar{P}^*$ of the closure.  The facets of $B(P)$ bijectively
correspond to the vertices of~$P$.  As an equation of $f$-vectors this reads as follows.
\begin{equation}\label{eq:bounded-f}
  f_k(\bounded{P}) \ = \ f_{n-k-1}(B(\bar{P}^*))-f_{n-k-1}(\partial B(\bar{P}^*)) \ = \ f_{n-k-1}(\interior{B(\bar{P}^*)})
\end{equation}
Moreover, since $h(\interior{B(\bar{P}^*)})$ is defined via the equation~\eqref{eq:h-from-f},
Proposition~\ref{prop:f-from-h} implies that
\begin{equation}\label{eq:h-h}
h_{n-k}(\bounded{P}) \ = \ h_k(\interior{B(\bar{P}^*)}) \ \whyrelation{\ref{prop:h-int}}{=} \ h_{n-k}(B(\bar{P}^*))\quad .
\end{equation}

\begin{exmp}
  A simplicial $n$-polytope is \emph{neighborly} if any set of $\lfloor n/2\rfloor$ vertices forms a face.  Examples are
  provided by the \emph{cyclic polytopes}, that is, the convex hulls of finitely many points on the \emph{moment curve}
  \[ t\mapsto(t,t^2,\dots,t^n) \quad .\] The definition of neighborliness readily implies that any triangulation of a
  neighborly simplicial polytope without additional vertices is asff.  Corollary~\ref{cor:defined-by-boundary-n-odd} now
  says that each triangulation of an even-dimensional neighborly simplicial polytope has the same $f$-vector.  Such
  polytopes are called \emph{equidecomposable}.
\end{exmp}

The next example will suitably be generalized in Section~\ref{sec:max}.

\begin{exmp}\label{exmp:octahedron}
  Any triangulation of a $3$-polytope without additional vertices is asff.  For instance, see the triangulation $\Theta$
  of the regular octahedron in Figure~\ref{fig:exmp}.  Here we have
  \[
  \begin{array}{lclclclclcl}
    f(\Theta)  &=& (6,13,12,4) \ ,  &\quad&  f(\partial\Theta) &=& (6,12,8) \ ,  &\quad& f(\interior{\Theta}) &=&
    (0,1,4,4) \ , \\
    h(\Theta)  &=& (1,2,1,0,0) \ ,  &\quad&  h(\partial\Theta) &=& (1,3,3,1) \ , &\quad& h(\interior{\Theta}) &=&
    (0,0,1,2,1) \ .
  \end{array}
  \]
\end{exmp}

\section{Tight Spans and Triangulations of Hypersimplices}

\noindent
A \emph{distance function} is a symmetric matrix with real coefficients and a zero diagonal.  We identify distance
functions with vectors in $\RR^{\binom{n}{2}}$ in a natural way.  A non-negative distance function $d$ is a
\emph{metric} if it satisfies the triangle inequality $d(i,k)\le d(i,j)+d(j,k)$. 

We recall some definitions from the introduction.  Each finite metric $d\in\RR^{\binom{n}{2}}$ gives rise to a pointed
unbounded polyhedron
\[ P_d \ = \ \SetOf{x\in\RR^n}{x_i+x_j\ge d(i,j)\text{ for all $i,j$}} \quad . \] The bounded subcomplex $T_d :=
\bounded{P_d}$ is called the \emph{tight span} of~$d$.  The metric $d$ is \emph{generic} if the polyhedron~$P_d$ is
simple.

The \emph{second hypersimplex} \[ \Delta_{n,2} \ = \ \conv\SetOf{e_i+e_j}{1\le i<j\le n} \] is an $(n-1)$-polytope which
is not simplicial.  In fact, its facets are either $(n-2)$-simplices or $(n-2)$-dimensional hypersimplices
$\Delta_{n-1,2}$.  As in De Loera, Sturmfels, and Thomas~\cite{MR1357285} we will use graph theory language in order to
describe a regular polyhedral subdivision $\Delta^d$ of $\Delta_{n,2}$ induced by the metric~$d$: If we identify the
vertices of $\Delta_{n,2}$ with the edges of the complete graph $K_n$ in a natural way then the cells of $\Delta^d$
correspond to subgraphs $\Gamma$ of $K_n$ (represented by their edge sets) which admit a height function
$\lambda\in\RR^n$ satisfying
\[
\lambda_i+\lambda_j=d(i,j) \text{ if $\{i,j\}$ is an edge}
\quad\text{and}\quad
\lambda_i+\lambda_j>d(i,j) \text{ if $\{i,j\}$ is not an edge of~$\Gamma$}\quad .
\]
The metric $d$ is generic if and only if $\Delta^d$ is a (regular) triangulation.  Conversely, each regular
triangulation of $\Delta_{n,2}$ gives rise to a generic metric. Hence in the generic case we can apply the results from
the previous sections.

In the next few steps we will explore the structure of $T_d$ in terms of the dual simplicial ball $\Delta^d$.  To this
end it is instrumental to begin with detailed information about the dual graph of~$\Delta_{n,2}$.  The small cases are,
of course, special: $\Delta_{3,2}$ is a triangle, and $\Delta_{4,2}$ is an octahedron, as studied in
Example~\ref{exmp:octahedron}.  The following is known, which is why we omit the (simple) proof.

\begin{lem}\label{lem:hypersimplex-dual-graph}
  Let $n\ge 5$.  Then the second hypersimplex $\Delta_{n,2}$ has $n$~facets isomorphic with $\Delta_{n-1,2}$ and $n$
  simplex facets.  Any two facets of hypersimplex type are adjacent, and their intersection is isomorphic with
  $\Delta_{n-2,2}$.  No two simplex facets are adjacent.  Each simplex facet is adjacent to $n-1$ hypersimplex facets.
\end{lem}

A consequence of this observation is that all the faces of a hypersimplex are either hypersimplices or simplices.

\begin{prop}\label{prop:inductive-step}
  For $n\ge 5$ let $\Delta$ be a triangulation of $\Delta_{n,2}$ such that on each $m$-dimensional hypersimplex face a
  triangulation with the same $f$-vector $(f^{(m)}_0,\dots,f^{(m)}_m)$ is induced.  Then we obtain
  \[
  f_{n-2}(\partial\Delta) \ = \ n + n f^{(n-2)}_{n-2}
  \quad \text{and} \quad
  f_k(\partial\Delta) \ = \ \sum_{i=1}^{n-1-k} (-1)^{i-1} \binom{n}{i} f^{(n-1-i)}_k \; \text{for $k<n-2$} \quad .
  \]
\end{prop}

\begin{proof}
  The claim for $f_{n-2}$ follows from the fact that $\Delta_{n,2}$ has $n$ simplex facets and $n$ hypersimplex facets,
  and that we assumed that each hypersimplex facet is triangulated into $f^{(n-2)}_{n-2}$ simplices of dimension~$n-2$.
  Lemma~\ref{lem:hypersimplex-dual-graph} says that the subgraph of the dual graph of $\Delta_{n,2}$ induced on the
  hypersimplex facets is a complete graph $K_n$.  Moreover, each face of dimension less than $n-2$ arises as a subface
  of a hypersimplex facet.  Therefore only the triangulations of the hypersimplex facets have to be taken into account,
  where doubles have to be removed.  The claim then follows from a standard inclusion-exclusion argument.
\end{proof}

Clearly, Proposition~\ref{prop:inductive-step} translates into various equations for the $g$- and $h$-vectors.  We
choose to establish the following relation.

\begin{cor}\label{cor:inductive-step}
   For $n\ge 5$ let $\Delta$ be a triangulation of $\Delta_{n,2}$ such that on each $m$-dimensional hypersimplex face a
   triangulation with the same $f$-vector $(f^{(m)}_0,\dots,f^{(m)}_m)$ is induced.  Then we obtain 
   \[
   g_k(\partial\Delta) \ = \ \sum_{i=1}^n\sum_{j=0}^{\min(i,k)}(-1)^{i+j-1}\binom{n}{i}\binom{i}{j}h^{(n-1-i)}_{k-j} \quad
   \text{for $k\le\lfloor n/2\rfloor$.}
   \]
   Here $(h^{(k)}_0,\dots,h^{(k)}_k)$ denotes the common $h$-vector of the $k$-faces.
 \end{cor}

\begin{proof}
  \begin{align*}
    g_k(\partial\Delta) \ & = \ \sum_{i=0}^k(-1)^{k-i}\binom{n-i}{k-i}f_{i-1}(\partial\Delta) \\
    & \whyrelation{\ref{prop:inductive-step}}{=} \
          \sum_{i=0}^k(-1)^{k-i}\binom{n-i}{k-i} \, \left(\sum_{j=1}^{n-i}(-1)^{j-1}\binom{n}{j}f^{(n-1-j)}_{i-1}\right) \\
    & = \ \sum_{i=0}^k(-1)^{k-i}\binom{n-i}{k-i} \, \left(\sum_{j=1}^{n}(-1)^{j-1}\binom{n}{j}f^{(n-1-j)}_{i-1}\right)
    \quad \text{(since $f^{(n-1-j)}_{i-1}=0$ if $j>n-i$)} \\
    & = \ \sum_{i=0}^k\sum_{j=1}^{n}\binom{n}{j}(-1)^{k-i+j-1}\binom{n-i}{k-i} f^{(n-1-j)}_{i-1}\\
    & = \ \sum_{i=0}^k\sum_{j=1}^{n}\binom{n}{j}(-1)^{k-i+j-1}\left(\sum_{l=0}^j\binom{j}{l}\binom{n-j-i}{k-l-i}\right)
    f^{(n-1-j)}_{i-1} \\
    & = \ \sum_{j=1}^{n} \sum_{l=0}^j (-1)^{j+l-1}\binom{n}{j}\binom{j}{l} \sum_{i=0}^k (-1)^{k-l-i}\binom{n-j-i}{k-l-i}
    f^{(n-1-j)}_{i-1} \\
    & = \ \sum_{j=1}^{n} \sum_{l=0}^j (-1)^{j+l-1}\binom{n}{j}\binom{j}{l} \sum_{i=0}^{k-l} (-1)^{k-l-i}\binom{n-j-i}{(n-j)-(k-l)}
    f^{(n-1-j)}_{i-1} \\
    & = \  \sum_{j=1}^{n} \sum_{l=0}^j (-1)^{j+l-1}\binom{n}{j}\binom{j}{l} h^{(n-1-j)}_{k-l} \\
    & = \  \sum_{j=1}^{n} \sum_{l=0}^{\min(j,k)} (-1)^{j+l-1}\binom{n}{j}\binom{j}{l} h^{(n-1-j)}_{k-l} \quad .
  \end{align*}
\end{proof}




We call a distance function $e\in\RR^{\binom{n}{2}}$ \emph{isolated} if there is an index $i\in\{1,\dots,n\}$ and a (not
necessarily positive) real number $\lambda\ne 0$ such that $e(i,j)=e(j,i)=\lambda$ for all $j\ne i$ and $e(j,k)=0$
otherwise.  Moreover, we say that two metrics are \emph{equivalent} if they differ by a linear combination of isolated
distance functions.  The following is known.

\begin{prop}\label{prop:equivalent-metrics}
  Let $d$ be a generic metric.
  \begin{enumerate}
  \item If $d$ and $d'$ are equivalent metrics then $\Delta^d=\Delta^{d'}$.
  \item For each generic metric $d$ there is a unique equivalent generic metric $d'$ such that $B(P_{d'})$ is
    combinatorially equivalent to $\Delta^{d'}=\Delta^d$.\label{it:ideal}
  \end{enumerate}
\end{prop}

An metric $d'$ is \emph{ideal} if it satisfies $\Delta^{d'}\cong B(P_{d'})$.
Proposition~\ref{prop:equivalent-metrics}\eqref{it:ideal} then reads as: Each generic metric is equivalent to an ideal
one.  The equivalence class of metrics of an ideal generic metric $d'$ on $n$ points can be described as follows: The
triangulation $\Delta^{d'}$ induces a triangulation of the boundary of the hypersimplex $\Delta_{n,2}$.  For $n\ge 5$,
$\Delta_{n,2}$ has $n$ simplex facets, and the simplicial balls $B(P_d)$ corresponding to non-ideal metrics equivalent
to $d'$ arise from $\Delta^{d'}=\Delta^d$ by gluing additional $(n-1)$-simplices to the simplex facets of
$\Delta_{n,2}$.

\begin{exmp}\label{exmp:4points}
  Consider the metric on four points given by the matrix
  \begin{equation}\label{eq:exmp}
    d=
    \begin{pmatrix}
      0 & 2 & 3 & 2\\
      2 & 0 & 2 & 3\\
      3 & 2 & 0 & 2\\
      2 & 3 & 2 & 0
    \end{pmatrix} \quad .
  \end{equation}
  The metric~$d$ turns out to be generic, and the tight span $T_d=\bounded{P_d}$ is $2$-dimensional.  The corresponding
  simplicial ball $\Delta^d$ is a triangulation of the regular octahedron, that is, the hypersimplex $\Delta(4,2)$.  See
  Figure~\ref{fig:exmp}.

  \begin{figure}[htbp]\centering
    \includegraphics[width=.2\textwidth]{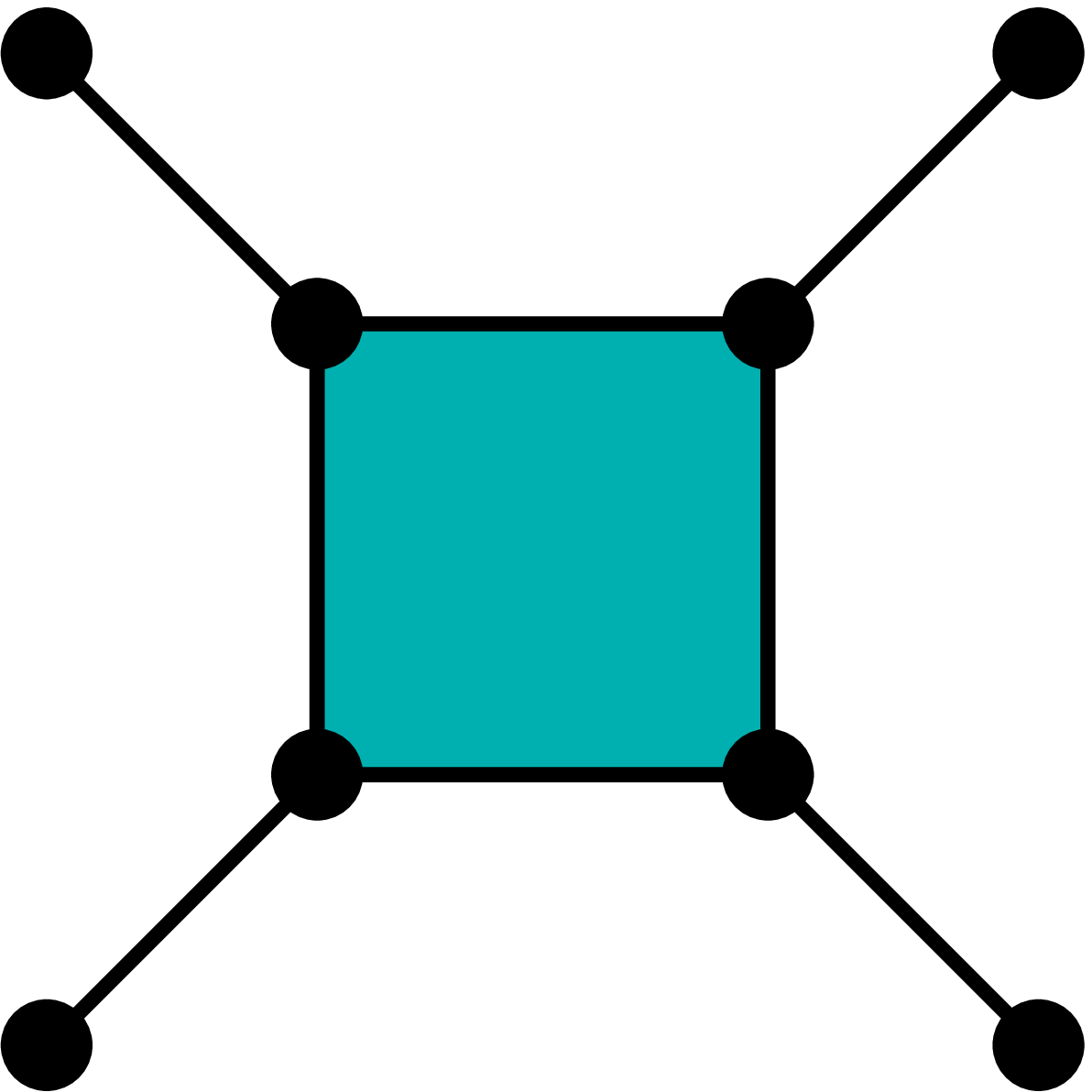}
    \qquad
    \includegraphics[width=.2\textwidth]{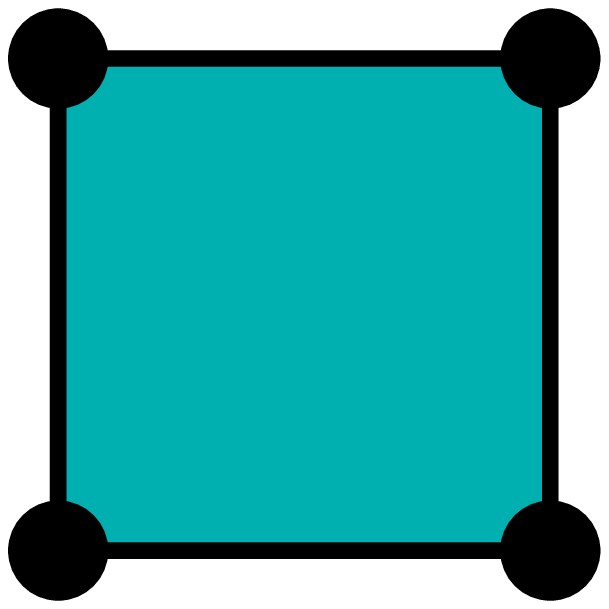}
    \qquad
    \includegraphics[width=.2\textwidth]{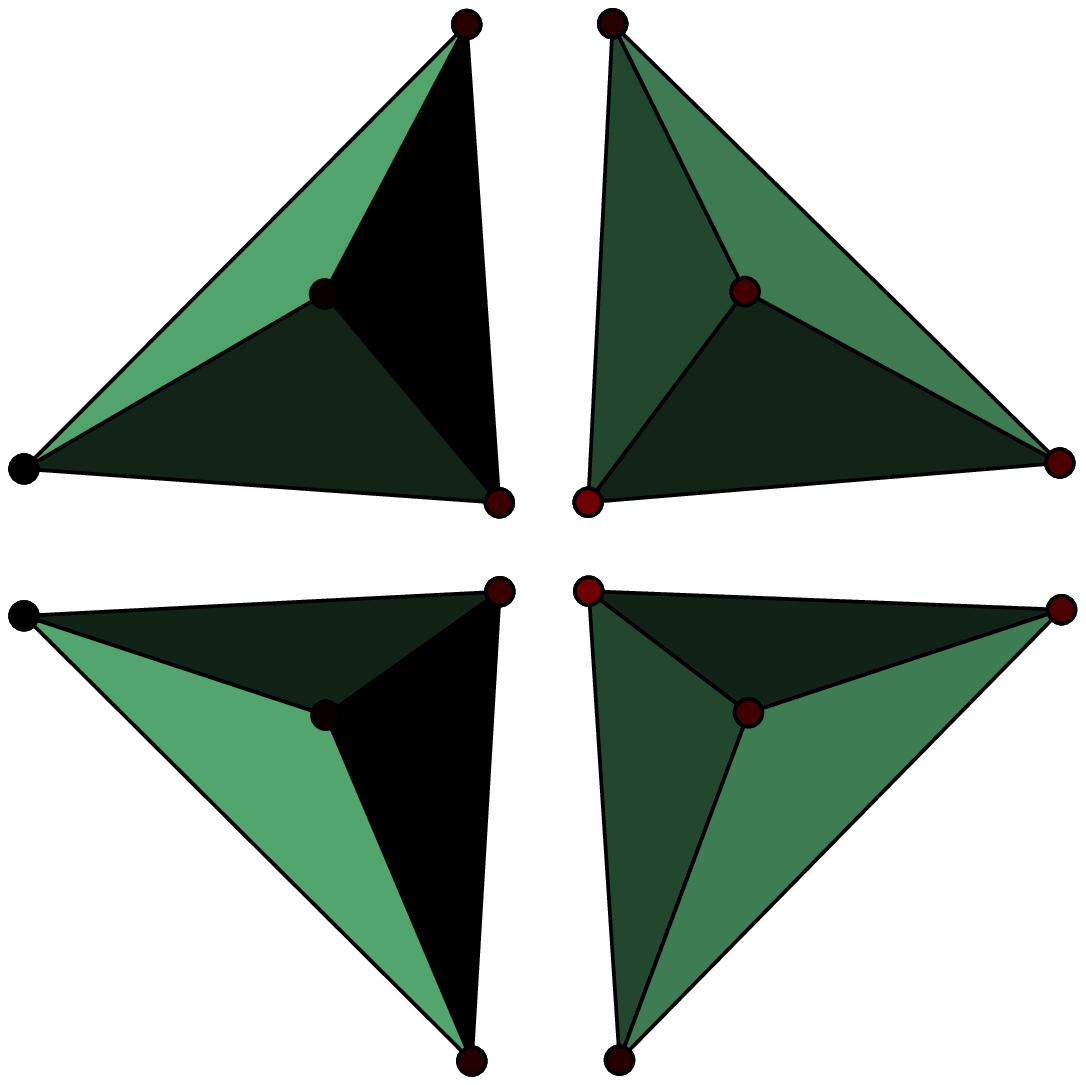}

    \caption{The tight span of the metric~$d$ defined in~\eqref{eq:exmp}, the tight span of an equivalent ideal metric
      $d'$, and the corresponding triangulation $\Delta^d=\Delta^{d'}$. Images produced with \texttt{polymake}~\cite{polymake} and
    \texttt{JavaView}~\cite{javaview}\label{fig:exmp}}
  \end{figure}

  The metric
  \[
  d'=
  \begin{pmatrix}
    0 & 1 & 2 & 1\\
    1 & 0 & 1 & 2\\
    2 & 1 & 0 & 1\\
    1 & 2 & 1 & 0
  \end{pmatrix}
  \]
  is equivalent to~$d$ and ideal, that is, its tight span satisfies $T_{d'}\cong\Delta^{d'}=\Delta^d$.
\end{exmp}

\begin{lem}\label{lem:ideal}
  Let $d,d'\in\RR^{\binom{n}{2}}$ be equivalent metrics such that $d'$ is ideal.  Then $h_k(T_d)=h_k(T_{d'})$ for
  $k\ne1$ and $h_1(T_d)\le h_1(T_{d'})+n$.
\end{lem}

Throughout the following we consider a fixed generic metric~$d$.

We summarize results of Develin~\cite{Develin}.  As before we identify a metric $d$ on $n$ points with an element of
$\RR^{\binom{n}{2}}$ and a graph on $n$ nodes with a $0/1$-vector of the same length $\binom{n}{2}$.

\begin{defn}
  For a given weight vector $w\in\RR_+^n$ on $n$ points we call a non-negative vector $\mu\in\RR^{\binom{n}{2}}$ a
  \emph{fractional $w$-matching} if $\sum_{i=1}^{n}\mu(i,j)=w_j$ for all $1\le j\le n$.  The \emph{support} $\supp\mu$ is
  the graph of those edges $(i,j)$ with $\mu(i,j)>0$.
\end{defn}

For a given graph $\Gamma\in\{0,1\}^{\binom{n}{2}}$, and $w\in\RR_+^n$ with $w_i=\deg_\Gamma(i)$, consider the linear
program
\begin{equation}\label{eq:LP}
  \begin{array}{ll}
    \max\ \langle \mu,d \rangle & \text{subject to}\\ [0.5ex]
    \sum_{i=1}^n \mu(i,j)=w_j   & \text{for all $1\le j\le n$, and}\\
    \mu\ge 0                    & \hspace{10em} .
  \end{array}
\end{equation}

A fractional $w$-matching is called \emph{optimal} if it is an optimal solution of this linear program.

\begin{thm}{\rm (Develin~\cite{Develin})}\label{thm:Develin}
  Let $d$ be a generic metric on $n$ points.
  \begin{enumerate}
  \item For each graph $\Gamma\in\{0,1\}^{\binom{n}{2}}$ the linear program~\eqref{eq:LP} has a unique optimal
    solution~$\opt{\mu}(\Gamma)$.\label{it:opt}
  \item The graphs $\Gamma$ with $\supp\opt{\mu}(\Gamma)=\Gamma$ are precisely the cells of~$\Delta^d$.\label{it:cells}
  \item A cell $\Gamma$ is an interior simplex if and only if it is a spanning subgraph of~$K_n$ which is not isomorphic
    with the star $K_{1,n-1}$. \label{it:spanning}
  \item The support of an optimal $w$-matching for an arbitrary  $w\in\RR_+^n$ is a cell of $\Delta^d$.\label{it:support}
  \item No cell $\Gamma$ contains an non-trivial even tour. \label{it:even-cycle}
  \item The dimension of $T_d$ is bounded by
    \[ \lceil n/3 \rceil \ \le \ \dim T_d \ \le \ \lfloor n/2 \rfloor \quad . \] \label{it:bounds}
  \end{enumerate}
\end{thm}

Here a \emph{tour} in the graph~$\Gamma$ is any closed path $(v_0,v_1,\dots,v_m=v_0)$; it is \emph{trivial} if each of
its edges occurs at least twice.  A \emph{cycle} is a tour in which each edge occurs only once.  In particular,
statement~\eqref{it:even-cycle} in the theorem implies that each vertex is contained in at most one cycle (which must
further be odd, if it exists).  Further, it turns out that the property~\eqref{it:even-cycle} characterizes the
non-degenericity of~$d$; see \cite[Proposition~2.10]{Develin}.

The following lemma is a key step in obtaining upper bounds on the $f$-vectors of tight spans.  It gives a bound on the
number of facets of~$T_d$ in the case where the dimension $\dim T_d=\lfloor n/2 \rfloor$ is maximal.

\begin{lem}\label{lem:n/2-bound}
  The triangulation $\Delta^d$ is asff.  Moreover,
  \[
  f_{\lfloor n/2\rfloor}(T_d) \ = \ f_{\lceil n/2\rceil-1}(\interior{\Delta^d})  \ \le \
  \begin{cases}
    1 & \text{if $n$ even}\\
    n & \text{if $n$ odd} \quad .
  \end{cases}
  \]
\end{lem}

\begin{proof}
  Any spanning subgraph of the complete graph $K_n$ needs at least $\lceil n/2 \rceil$ edges.  In view of
  Theorem~\ref{thm:Develin}\eqref{it:spanning} this implies that an interior face of $\Delta^d$ is at least of dimension
  $\lceil n/2 \rceil-1=\lfloor (n-1)/2\rfloor$ or, equivalently, that $\Delta^d$ is asff.
  
  Assume first that $n$ is even, and that $\Gamma$ is a graph with $n/2$ edges which corresponds to an interior simplex
  of~$\Delta^d$.  This says that $\Gamma$ is a perfect matching of~$K_n$ and hence an optimal solution of the linear
  program~\eqref{eq:LP} for the weight $w=(1,1,\dots,1)$.  From the uniqueness result
  Theorem~\ref{thm:Develin}\eqref{it:opt} it thus follows that $f_{n/2-1}(\interior{\Delta^d})\le 1$.
  
  Now let $n$ be odd.  Then $\Gamma$ is a spanning subgraph of~$K_n$ with $(n+1)/2$ edges.  This implies that $\Gamma$
  has a unique node $t$ of degree~$2$.  Clearly, there are $n$ choices for $t$.
\end{proof}

Note that $\Delta^d$ being asff is equivalent to the upper bound $\dim T_d \le \lfloor n/2 \rfloor$ in
Theorem~\ref{thm:Develin}\eqref{it:bounds}.

As a further piece of notation we introduce
\begin{align*}
  H_k(n) \
   :=& \ \max\SetOf{h_k(\Delta)}{\text{$\Delta$ regular triangulation of $\Delta_{n,2}$}} \\
    \whyrelation{\eqref{eq:h-h}}{=}& \ \max\SetOf{h_k(T_d)}{\text{$d$ ideal metric on $n$ points}} \quad .
\end{align*}
We are now ready to prove our main result.

\begin{thm}\label{thm:upper-bound}
  The $h$-vector of a regular triangulation $\Delta$ of the hypersimplex $\Delta_{n,2}$ is bounded from above by
  \[
  H_k(n) \ \le \ \binom{n}{2k} \quad \text{for $k\ne 1$}
  \]
  and $H_1(n) \le \binom{n}{2}-n$.
\end{thm}

Via Proposition~\ref{prop:f-from-h} this upper bound on the $h$-vector gives the recursion
\[
F_k(n) \ = \ 2F_k(n-1)+F_{k-1}(n-2) \quad ,
\]
where $F_k(n)$ is the maximal number of $k$-faces of the tight span of any generic metric on $n$ points.  This further
translates into the following equivalent upper bound for the $f$-vector:
\[ F_k(n) \ \le \ 2^{n-2k-1}\frac{n}{n-k}\binom{n-k}{k} \quad . \] In Section~\ref{sec:max} it will be shown that these
bounds are tight.  There even is a regular triangulation of~$\Delta_{n,2}$ which simultaneously maximizes all entries of
the $h$-vector.  Note that this fact will be used in the proof of this theorem.

The bound $F_0(n) \le 2^{n-1}$ for the number of vertices of a tight span also follows from the known fact that the
normalized volume of $\Delta_{n,2}$ equals $2^{n-1}-n$: The vertices of a tight span of an ideal generic metric are in
$1-1$ correspondence with the facets of a regular triangulation of~$\Delta_{n,2}$; and changing from the ideal metric to
an equivalent non-ideal metric allows for another $n$ vertices in the tight span.  As there are unimodular (and regular)
triangulations of $\Delta_{n,2}$, for instance, the thrackle triangulations studied by De Loera, Sturmfels, and
Thomas~\cite{MR1357285}, it is clear that this bound is tight.

We need some elementary facts about multinomial coefficients, which we phrase as equations of binomial coefficients.
Moreover, it will be convenient to make use of \emph{Kronecker's delta} notation
\[
\delta_{i,k} \ = \ \begin{cases} 1 & \text{if $i=k$,} \\ 0 & \text{otherwise.} \end{cases}
\]

\begin{lem}\label{lem:main:a}
  \[
  \sum_{i=1}^n (-1)^{i+k} \binom{n}{i} \binom{i}{k-1} (n-i) \ = \ n\,\delta_{1,k}
  \]
\end{lem}

\begin{proof}
  For $k=0$ we have $\binom{i}{-1}=0$, and the claim is obvious.  So we assume that $k>0$.
  \begin{align*}
    \sum_{i=1}^n (-1)^{i+k} \binom{n}{i} \binom{i}{k-1} (n-i) \ &= \
    \sum_{i=k-1}^n (-1)^{i+k} \binom{n}{i} \binom{i}{k-1} (n-i) \; - \; \delta_{1,k}(-1)^1 \binom{n}{0} \binom{0}{0} \\
    &= \ k\binom{n}{k} \sum_{i=k-1}^n (-1)^{i+k} \binom{n-k}{i-(k-1)} \; + \; n\delta_{1,k} \\
    &= - \ k\binom{n}{k} \sum_{i=0}^{n-(k-1)} (-1)^k \binom{n-k}{i} \; + \; n\delta_{1,k} \\
    &= \ n\delta_{1,k} \quad .
  \end{align*}
\end{proof}

\begin{lem}\label{lem:main:b}
  \[ \sum_{i=j}^n (-1)^{i+j-1}\binom{n}{i}\binom{i}{j}\binom{n-i}{2(k-j)} \ = \ 0 \quad . \]
\end{lem}

\begin{proof}
  \begin{align*}
    \sum_{i=j}^n (-1)^{i+j-1}\binom{n}{i}\binom{i}{j}\binom{n-i}{2(k-j)} \
    =& \ \binom{n}{j}\binom{n-j}{2(k-j)} \sum_{i=j}^n (-1)^{i+j-1}\binom{n-2k+j}{i-j} \\
    =& \ - \binom{n}{j}\binom{n-j}{2(k-j)} \sum_{i=0}^{n-j} (-1)^k \binom{n-2k+j}{i} \\
    =& \ - \binom{n}{j}\binom{n-j}{2(k-j)} \sum_{i=0}^{n-2k+j} (-1)^k \binom{n-2k+j}{i} \\
    =& \ 0 \quad .
  \end{align*}
\end{proof}

\begin{proof}[Proof of Theorem~\ref{thm:upper-bound}]
  The hypersimplex $\Delta(4,2)$ is the regular octahedron, and (up to combinatorial equivalence) it has a unique
  triangulation $\Theta$ without additional vertices; see the Examples~\ref{exmp:octahedron} and~\ref{exmp:4points}.
  Then $h(\Theta)=(1,2,1,0,0)$.  This settles the case $n=4$.

  We will proceed by induction on $n$.  From Proposition~\ref{prop:asff-h-vanishing} and Equation~\eqref{eq:h-h} it
  follows that maximizing the $h$-vector of $\Delta$ amounts to the same as maximizing the $g$-vector of the boundary
  $\partial\Delta$.  Hence, inductively we can assume that each hypersimplex $l$-face of $\Delta_{n,2}$ is maximally
  triangulated, that is, in the notation of Corollary~\ref{cor:inductive-step},
  $h^{(l)}_k=\binom{l+1}{2k}-(l+1)\delta_{1,k}$ for all~$k$.  We can write this as an equation rather than an inequality
  since we know from the construction in Section~\ref{sec:max} that this bound is attained.

  \begin{align*}
    h_k(\Delta)
    \ \whyrelation{\ref{prop:asff-h-vanishing}}{=}& \ g_k(\partial\Delta) \\
    \ \whyrelation{\ref{cor:inductive-step}}{=}& \
    \sum_{i=1}^n\sum_{j=0}^{\min(i,k)}(-1)^{i+j-1}\binom{n}{i}\binom{i}{j}h^{(n-1-i)}_{k-j} \\
    \ =& \
    \sum_{i=1}^n\sum_{j=0}^{\min(i,k)}(-1)^{i+j-1}\binom{n}{i}\binom{i}{j}\left[\binom{n-i}{2(k-j)}-(n-i)\delta_{1,k-j}\right] \\
     =& \ \sum_{i=1}^n\sum_{j=0}^{\min(i,k)}(-1)^{i+j-1}\binom{n}{i}\binom{i}{j}\binom{n-i}{2(k-j)} -
     \sum_{i=1}^n\sum_{j=0}^{\min(i,k)}(-1)^{i+j-1}\binom{n}{i}\binom{i}{j}(n-i)\delta_{1,k-j} \\
     =& \ \sum_{i=1}^n\sum_{j=0}^{\min(i,k)}(-1)^{i+j-1}\binom{n}{i}\binom{i}{j}\binom{n-i}{2(k-j)}
     \; - \; \sum_{i=1}^n (-1)^{i+k}\binom{n}{i}\binom{i}{k-1}(n-i)\\
     \whyrelation{\ref{lem:main:a}}{=}& \ \sum_{j=0}^{k}\sum_{i=1}^n (-1)^{i+j-1}\binom{n}{i}\binom{i}{j}\binom{n-i}{2(k-j)} \; - \; n \delta_{1,k} \\
     =& \ \sum_{j=0}^{k}\sum_{i=j}^n (-1)^{i+j-1}\binom{n}{i}\binom{i}{j}\binom{n-i}{2(k-j)}
     \; - \; (-1)^{-1}\binom{n}{0}\binom{0}{0}\binom{n-0}{2(k-0)} \; - \; n \delta_{1,k} \\
     \whyrelation{\ref{lem:main:b}}{=}& \ \binom{n}{2k} \; - \; n \delta_{1,k}
  \end{align*}
\end{proof}

\section{A Metric with Maximal $f$-Vector}
\label{sec:max}

\noindent
In the sequel we will prove that the upper bounds given are tight.  To this end, for each $n\ge 4$, we define the metric
$\dmax^n$ by letting
\[
\dmax^n(i,j) \ = \ 1 + \frac{1}{n^2+in+j} \quad ,
\]
for $1\le i < j \le n$.  We suitably abbreviate $\Pmax^n=P_{\dmax^n}$ and $\Tmax^n=T_{\dmax^n}$.

\begin{prop}
  The metric $\dmax^n$ is generic.
\end{prop}

\begin{proof}
  Due to \cite[Proposition~2.10]{Develin} it suffices to show that no graph $\Gamma$ corresponding to a cell of
  $\Delta^d$ contains a non-trivial even tour.  Assuming the contrary, let $C=(i_1,i_2,\dots,i_{2n},i_1)$ be such a
  tour.  Then we have a non-trivial affine dependence
  \[
  \sum_{(k,l) \in A}\dmax^n(k,l)= \sum_{(k,l) \in B}\dmax^n(k,l)
  \]
  with $A=\{(i_1,i_2),(i_3,i_4),\dots\}$ and $B=\{(i_2,i_3),\dots,(i_{2n-2},i_{2n-1}),(i_{2n},i_1)\}$.  But this
  contradicts the fact that $\{\dmax^n(i,j)\}$ is a linearly independent set over $\QQ$.
\end{proof}

\begin{figure}[htbp]\centering
  \includegraphics[height=.35\textwidth]{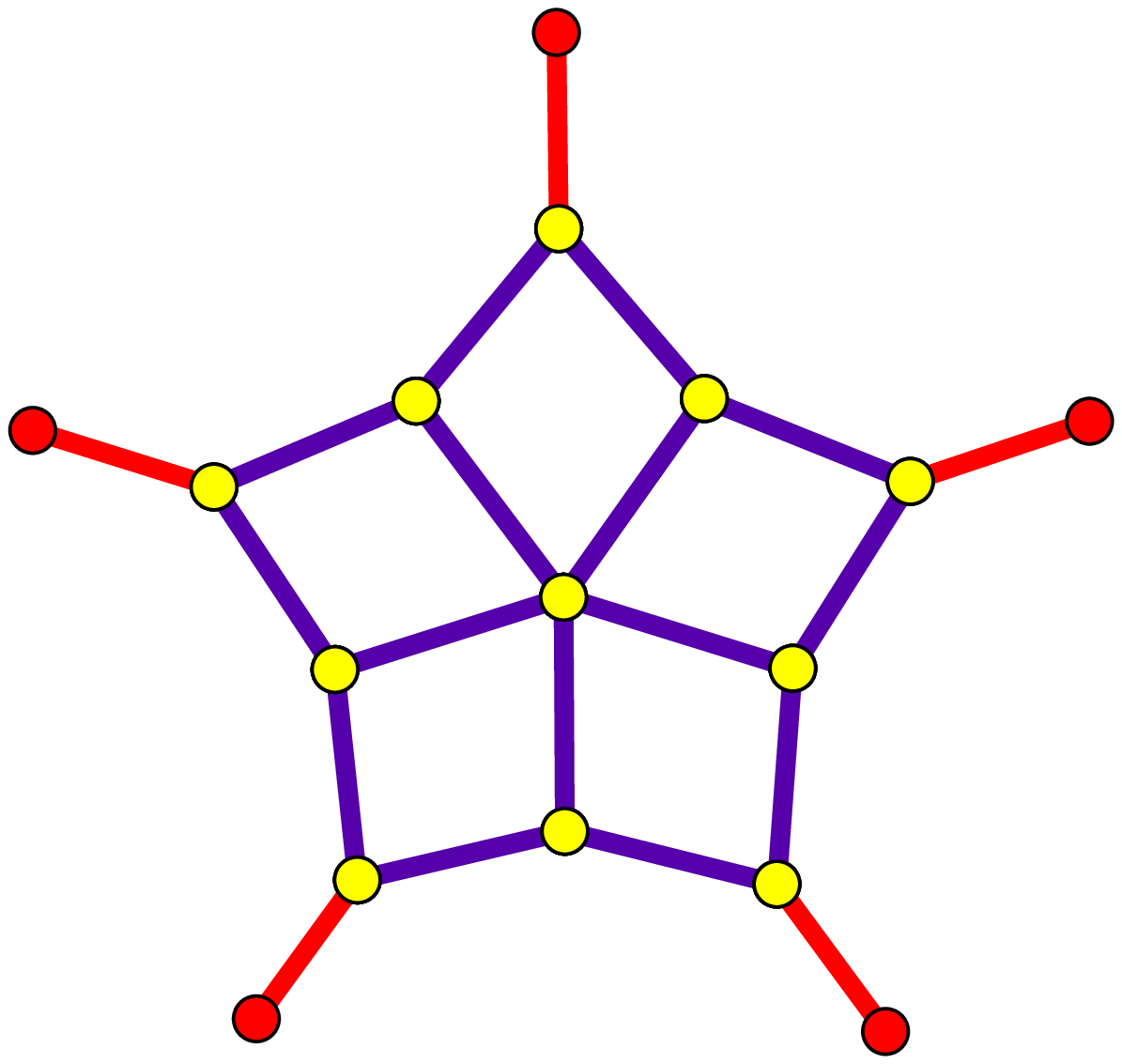}\quad
  \includegraphics[height=.35\textwidth]{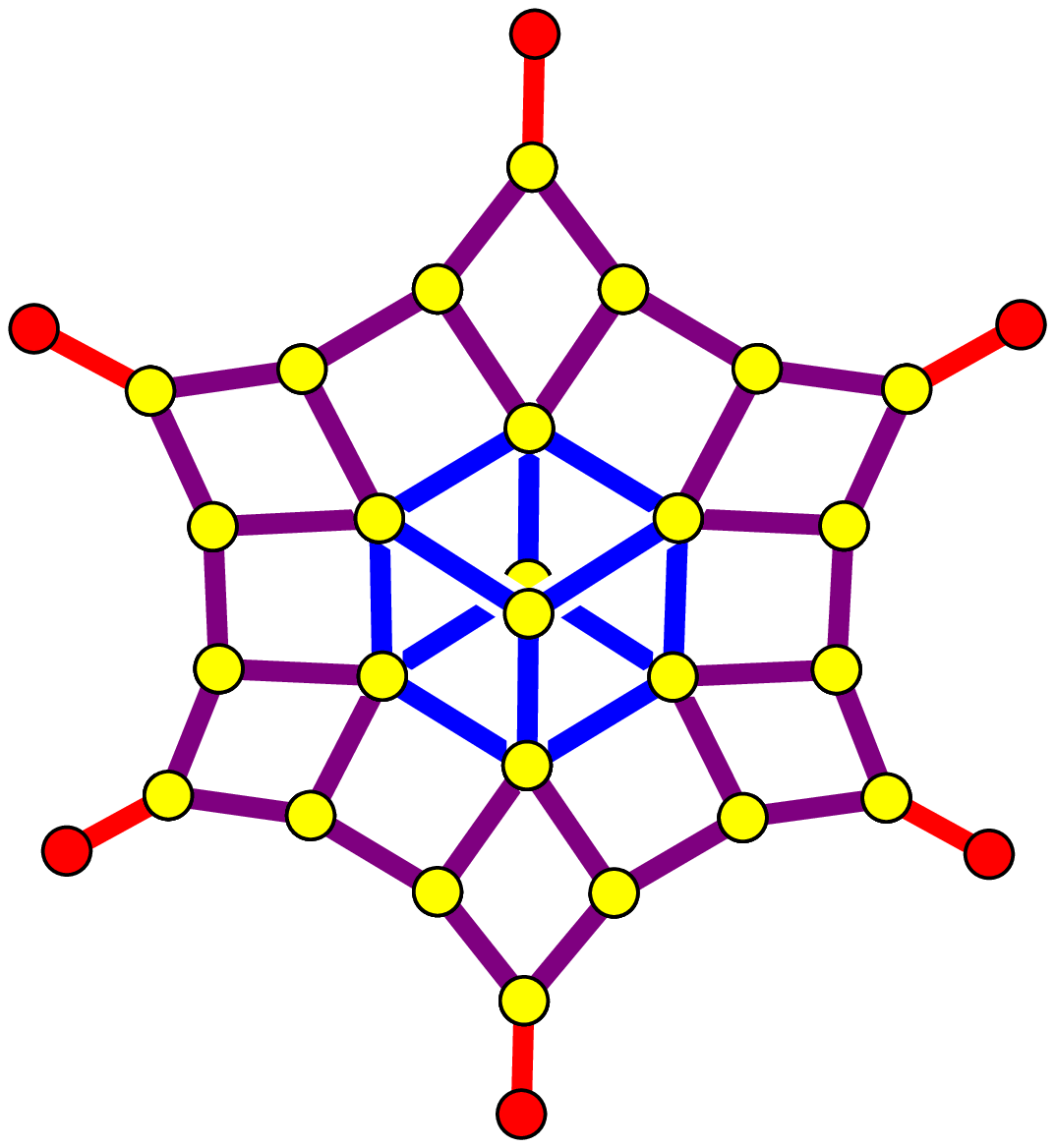}

  \caption{Visualization of (the graphs of) the tight spans $\Tmax^5$, with $f$-vector $(16,20,5)$, and $\Tmax^6$, with
    $f$-vector $(32,48,18,1)$.  The unique $3$-face of $\Tmax^6$ is a cube.  The two corresponding triangulations occur
    under the name ``thrackle triangulations'' in De Loera, Sturmfels, and Thomas~\cite{MR1357285}.  Moreover,
    $\Tmax^6$, or rather the tight span of an equivalent ideal metric, is $\#66$ in Sturmfels and Yu~\cite{MR2097310}.}
\end{figure}

The key property of the metric $\dmax^n$ is the following.

\begin{lem}\label{lem:d-max-property}
  For $1\leq i\leq j\leq k \leq l\leq n$ we have
  \[
  \dmax^n(i,j)-\dmax^n(i,k) \ \leq \ \dmax^n(j,l)-\dmax^n(k,l)\phantom{\quad.}
  \]
  and
  \[
  \dmax^n(i,l)-\dmax^n(i,k) \ \leq \ \dmax^n(j,l)-\dmax^n(j,k)\quad.
  \]
\end{lem}

\begin{proof}
  Without loss of generality we can assume $i<j<k<l$.  Then we have
  \begin{align*}
    \dmax^n(i,j)-\dmax^n(i,k) \ =& \ \frac{1}{n^2+in+j} - \frac{1}{n^2+in+k}\\
    =& \ \frac{k-j}{(n^2+in+j)(n^2+in+k)} \ < \ \frac{(k-j)n}{(n^2+jn+l)(n^2+kn+l)}\\
    =& \  \frac{1}{n^2+jn+l} - \frac{1}{n^2+kn+l}
    \ = \ \dmax^n(j,l)-\dmax^n(k,l) \quad .
  \end{align*}
  The other inequality follows from a similar computation.
\end{proof}

It is clear that also all submetrics of $\dmax^n$, that is, metrics induced on subsets of $\{1,\dots,n\}$, share this
property.  To further analyze $\dmax^n$ and its tight span we require an additional characterization of the cells in the
tight span of a generic metric.  In the sequel we write $E(\Gamma)$ for the set of edges of a graph~$\Gamma$.
  
\begin{prop}\label{prop:connected-criteria}
  Let $d$ be a generic metric on $n$ points, and let $\Gamma$ be a connected graph with $n$ vertices, $n$ edges and without
  non-trivial even tours. Then $\Gamma$ defines a cell of $\Delta^d$ if and only if for all $\{v,w\} \not\in E(\Gamma)$ we have
  \begin{equation}\label{eq:fo}
  d(v,w) \ \leq \ \sum_{k=1}^{m-1} (-1)^{k-1} d(v_k,v_{k+1}),
  \end{equation}
  where $P=(v=v_1,v_2,\dots,v_m=w)$ is any path from $v$ to $w$ of odd length.
\end{prop}

\begin{proof}
  A connected graph with $n$ nodes and $n-1$ edges is a tree.  Therefore, $\Gamma$ can be seen as a tree with an
  additional edge which is contained in the unique (odd) cycle.  This implies that there is a path of odd length between
  any two vertices $v$ and $w$ (go around the odd cycle once if necessary).  While this path of odd length is not unique
  two such paths only differ by the insertion/deletion of trivial even tours or the direction in which the odd cycle is
  traversed.  Moreover, the set $P'$ of those edges occurring an odd number of times in the path~$P$ is independent of
  the choice of the path~$P$. A direct computation then shows that the value $\sum_{k=1}^{m-1} (-1)^{k-1}
  d(v_k,v_{k+1})$ is also independent of the choice of~$P$.
  
  Let $\Gamma$ be a cell of $\Delta^d$, and let $\{v,w\} \notin E(\Gamma)$ be an non-edge.  We consider the graph $C$
  consisting of $\{w,v\}$ and the edge set $P'$ of those edges which occur in the path $P$ an odd number of times.
  Clearly, $C$ is an even cycle in the complete graph, and we define $c'\in \RR^{\binom n2}$ as
  \begin{align}\label{cc}
    c'_{\alpha \beta} \ := \
    \begin{cases}
      1  &\text{for } \{\alpha,\beta\}=\{v_k,v_{k+1}\}\in E(C)\text{ and } k \text{ odd} \\
      -1 &\text{for } \{\alpha,\beta\}=\{v_k,v_{k+1}\}\in E(C)\text{ and } k \text{ even}\\
      1  &\text{for } \{\alpha,\beta\}=\{v,w\} \\
      0  &\text{otherwise}
    \end{cases}\quad.
  \end{align}
  Then $c:=\Gamma+\frac 1 2 c'$ is a feasible point of \eqref{eq:LP} and we have
  \[
  \langle c, d\rangle
  \ = \ \langle \Gamma, d \rangle+ \frac{1}{2} \left(-\sum_{k=1}^{m-1} (-1)^{k-1} d(v_k,v_{k+1}) +d(v,w)\right)
  \ > \ \langle\Gamma,d\rangle\quad.
  \]
  Since Theorem \ref{thm:Develin}\eqref{it:cells} establishes the optimality of~$\Gamma$ we can infer that the non-edge
  $\{v,w\}$ satisfies the inequality~\eqref{eq:fo}.
  
  For the reverse direction let $\Gamma$ be a graph such that \eqref{eq:fo} is true for all $\{v,w\} \notin E(\Gamma)$.
  Further let $\opt{\mu}(\Gamma)$ be the optimal solution to the linear program~\eqref{eq:LP}, which is unique due to
  Theorem \ref{thm:Develin}\eqref{it:opt}.  Then Theorem \ref{thm:Develin}\eqref{it:cells} tells us that we have to show
  $\opt{\mu}(\Gamma)=\Gamma$.  Assuming the converse, Theorem \ref{thm:Develin}\eqref{it:support} gives us $\{v,w\}\notin
  E(\Gamma)$ with $\opt{\mu}(\Gamma)_{vw}>0$. Let $c=\opt{\mu}(\Gamma)-\frac{\opt{\mu}(\Gamma)_{vw}}{2} c'$ with $C$ and
  $c'$ as in the first part of the proof. Then we have
  \[
    \langle c',d \rangle
    \ = \ \langle c,d \rangle + \frac{\opt{\mu}(\Gamma)_{vw}} 2 \left(\sum_{k=1}^{m-1} (-1)^{k-1} d(v_k,v_{k+1}) -
      d(v,w)\right)
    \ \geq \ 0 \quad .
  \]
  But this is a contradiction to the fact that $\opt{\mu}(\Gamma)$ is the unique optimal value.
\end{proof}

As mentioned previously, Lemma~\ref{lem:d-max-property} is the only property of $\dmax^n$ which actually matters.

\begin{lem}\label{lem:cycle-opt}
  Let $d$ be any generic metric on $n$ points for which the inequalities in Lemma \ref{lem:d-max-property} hold, for
  example, $d=\dmax^n$.  Then the cycle
  \[ C \ = \ (1,(n+1)/{2}+1,2,(n+1)/{2}+2,\dots,(n-1)/{2},n,(n+1)/2,1)\]
  is a cell of $\Delta^d$ if $n$ is odd.  If $n$ is even then the graph $D$ consisting of the cycle
  \[ C' \ = \ (1,n/{2}+2,2,n/{2}+3,\dots,n/{2}-1,n,n/{2},1) \]
  and the additional edge $\{1,n/{2}+1\}$ defines a cell.
\end{lem}

\begin{proof}
  We consider the case where $n$ is odd.  For each non-edge $\{j,l\}\notin E(C)$ we verify the conditions of Proposition
  \ref{prop:connected-criteria}.  The proof distinguishes four cases, the first of them being $j<l<(n+1)/2$.  The
  distance of $j$ and $l$ in the cycle~$C$ is even then, and as a path
  of odd length we can take
  \begin{align*}
    P \ = \ (&l,(n+1)/2+l,l+1,(n+2)/2+l+1,\dots,(n-1)/2,n,(n+1)/2,\\
             &1,(n+1)/2+1,2,(n+1)/2+2,\dots,j) \quad .
  \end{align*}
  Hence we have to show that        
  \begin{align}\label{eq:cycle-opt}
    d(j,l) \ \leq \ & \sum_{k=l}^{(n-1)/2}\left( d(k,(n+1)/2+k)-d(k+1,(n+1)/2+k)\right) \notag\\
                    & +\,d(1,(n+1)/2)\\
                    & - \, \sum_{k=1}^{j-1}\left(d(k,(n+1)/2+k)-d(k+1,(n+1)/2+k)\right) \quad \notag.
  \end{align}
  We compute
  \begin{align*}
    d(j,l)-d(1,(n+1)/2) \ = \ &\sum_{k=l}^{(n-1)/2} \left(d(j,k)-d(j,k+1)\right)\\
                          &- \, \sum_{k=1}^{j-1}\left( d(k,(n+1)/2)-d(k+1,(n+1)/2)\right)\quad.
  \end{align*}
  Considering the summands of the first sum, the first part of Lemma \ref{lem:d-max-property} yields
  \[
  d(j,k)-d(j,k+1) \ \leq \ d(k,(n+1)/2+k)-d(k+1,(n+1)/2+k)
  \]
  because $j\leq k\leq k+1\leq (n+1)/2 +k$.  The summands of the second sum satisfy $k\leq k+1\leq (n+1)/2\leq   
  (n+1)/2 +k$, whence the second part of Lemma \ref{lem:d-max-property} says that
  \[
  d(k,(n+1)/2)-d(k+1,(n+1)/2) \ \geq \ d(k,(n+1)/2+k)-d(k+1,(n+1)/2+k)\quad.
  \]
  By summing up we obtain the inequality~\eqref{eq:cycle-opt} as desired.
  
  The remaining three cases are $(n-1)/2<j<l$, $j<l-(n+1)/2<(n+1)/2+j<l$, and $l-(n-1)/2<j<(n+1)/2,(n-1)/2<l<j+(n+1)/2$.
  These, as well as the situation for $n$ even, are reduced to similar computations.
\end{proof}

\begin{figure}[htbp]\centering
  \begin{overpic}[width=.8\textwidth]{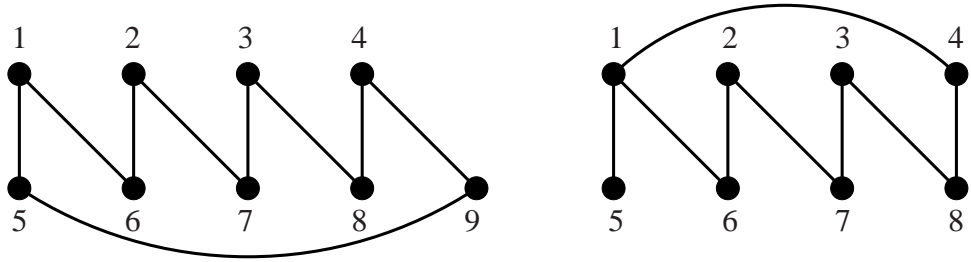}
    \put( 0.5 ,22){$1$}
    \put(12.25,22){$2$}
    \put(24   ,22){$3$}
    \put(35.75,22){$4$}

    \put( 0.5 ,3){$5$}
    \put(12.25,3){$6$}
    \put(24   ,3){$7$}
    \put(35.75,3){$8$}
    \put(47.5 ,3){$9$}

    \put(62.5 ,22){$1$}
    \put(74.25,22){$2$}
    \put(86   ,22){$3$}
    \put(97.75,22){$4$}

    \put(62.5 ,3){$5$}
    \put(74.25,3){$6$}
    \put(86   ,3){$7$}
    \put(97.75,3){$8$}
  \end{overpic}
  
  \caption{This illustrates Lemma~\ref{lem:cycle-opt}: Cycle~$C$ for $n=9$ odd (left) and graph $D$ for $n=8$ even (right).}
\end{figure}

\begin{thm}
  We have
  \[ h_i(\Tmax^n) \ = \ \binom{n}{2i} \quad . \]
\end{thm}

\begin{proof}
  First we show that we have equality in the bound of Lemma \ref{lem:n/2-bound} for $\dmax^n$ and all its submetrics. This is
  immediate from Lemma \ref{lem:cycle-opt} because for $n$ even the graph $D$ has a spanning subgraph with 
  $n/2$ edges corresponding to an interior simplex of $\Delta^d$ by Theorem \ref{thm:Develin}\eqref{it:spanning}. For $n$ odd we
  find $n$ spanning subgraphs of $C$ with $(n+1)/2$ edges each. These are exactly the bounds of Lemma \ref{lem:n/2-bound}. 
  
  Now the result follows from the computation in the proof of Theorem \ref{thm:upper-bound}.  
\end{proof}

\section{Towards a Lower Bound}
\noindent
Before we can prove something about lower bounds we require an additional lemma on the graphs defining cells of
$\Delta^d$, that is, graphs supporting optimal solutions of the linear program~\eqref{eq:LP}.

\begin{lem}\label{lem:B11}
  Let $w=(b,1,\dots,1)\in \RR^n$, and let $\Gamma$ be the support of the corresponding solution of the optimal
  fractional $w$-matching.  Then for the connected component $C$ of $\Gamma$ containing vertex $1$ exactly one of the
  following is true:
  \begin{enumerate}
  \item Either the component $C$ consists of one odd cycle and $b-1$ additional edges incident with the vertex~$1$,
    \label{it:B11-cycle}
  \item or the component $C$ consists of $b$ edges incident with the vertex~$1$.\label{it:B11-non_cycle}
  \end{enumerate}
  All other connected components of $\Gamma$ are isolated edges or odd cycles.
\end{lem}

\begin{proof}
  Let $\mu$ be a fractional $w$-matching with support graph $\Gamma$.  First, no vertex other than~$1$ can have a degree
  greater than or equal to~$3$: Suppose otherwise that there is a vertex~$x\ne 1$ with three neighbors $u,v,w$.  Since
  the total weight of the edges through~$x$ equals one, we have $\mu(x,u),\mu(x,v),\mu(x,w)<1$.  This implies that each
  of $u,v,w$ must be adjacent to another vertex (via an edge of weight less than one), and these paths continue further
  into all three directions starting from $x$.  Because the graph~$\Gamma$ is finite eventually these three paths must
  reach a vertex that they already saw previously.  Since we started into three directions it is not possible that all
  the vertices that we saw lie on one cycle.  Therefore there are at least two cycles in the connected component of~$x$,
  which implies that there is a non-trivial even tour through~$x$.  And this is forbidden by
  Theorem~\ref{thm:Develin}\eqref{it:even-cycle}.

  The same argument also shows that the vertex~$1$ is contained in at most one (odd) cycle.  Moreover, each vertex
  adjacent to~$1$ which is not contained in the odd cycle through~$1$ (if it exists) cannot be adjacent to any other
  vertex: Otherwise it would also generate a path which must end in a cycle as above.  Note that all edges in a cycle
  necessarily have weight $1/2$.

  If $\mu(x,y)=1$ for some $x,y\ne 1$ then both, $x$ and $y$ are only contained in the edge $\{x,y\}$.  Therefore the
  claim.
\end{proof}

\begin{figure}[htbp]\centering
  \includegraphics[width=.8\textwidth]{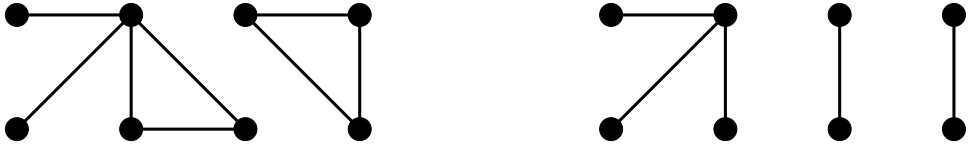}

  \caption{Graphs supporting an optimal $(b,1,\dots,1)$-matching as in Lemma~\ref{lem:B11} for $n=8$ and $b=4$ (two
    components to the left) and $b=3$ (three components to the right), respectively.}
\end{figure}

The case $b=1$ in the preceding result (with the same kind of argument) occurs in the proof of
Theorem~\ref{thm:Develin}\eqref{it:bounds} which is \cite[Theorem~3.1]{Develin}.

As a lower bound analog to Lemma \ref{lem:n/2-bound} for generic metrics we show the following theorem.  The three
different cases correspond to the congruence class of $\NN$ modulo~$3$.

\begin{thm}\label{thm:lower-bound}
  Let $d$ be a generic metric on $n$ points such that $T_d$ has dimension $\lceil n/3 \rceil$. Then we have
  \[
  f_{\lceil n/3 \rceil}(T_d) \ = \ f_{\lfloor 2n/3\rfloor-1}(\interior{\Delta^d})\ \geq \
  \begin{cases} 
    n\cdot3^{k-2}+3^k & \text{if $n=3k$} \\ 
    3^{k-1}           & \text{if $n=3k+1$} \\ 
    5\cdot 3^{k-1}    & \text{if $n=3k+2$}\quad.
  \end{cases}
  \]
\end{thm}

\begin{proof}
  Let first $n=3k+1$ and $\Gamma$ be the support of the optimal fractional $w$-matching for $w=(1,\dots,1)$.
  Lemma~\ref{lem:B11} yields that $\Gamma$ only consists of isolated edges and odd cycles. As $\Gamma$ cannot have a
  spanning subgraph with more than $\lfloor 2n/3\rfloor$ edges (since we assumed that $\dim T_d=\lceil n/3\rceil$) the
  only possibility is that $\Gamma$ consists of $k-1$ cycles of length three and two isolated edges. Since each
  $3$-cycle has exactly three spanning subgraphs we get at least $3^{k-1}$ faces of dimension $k=(n-1)/3$, as desired.
  
  For $n=3k$ a similar argument yields $3^k$ faces of dimension~$k$.  Additionally, we consider the support $\Gamma'$ of
  the optimal fractional $w$-matching for $w=(3,1,\dots,1)$, and again we can apply Lemma~\ref{lem:B11}.  If we were in
  case \eqref{it:B11-cycle} then $\Gamma'$ had a spanning subgraph of at most $3$ (from the connected component
  containing vertex $1$) plus $2(k-2)$ (from $k-3$ cycles of length three and two isolated edges in the rest) edges,
  summing up to $2k-1<2n/3$ altogether, which is impossible.  So we are in case \eqref{it:B11-non_cycle} of
  Lemma~\ref{lem:B11}.  Then we get spanning subgraphs of $\Gamma'$ with $3$ (connected component containing vertex $1$)
  plus $2k-3$ edges, which makes $2n/3$ altogether.  Again each of the $k-2$ cycles of length three of $\Gamma'$ has
  three possible spanning subgraphs yielding $3^{k-2}$ faces. These are all different from those obtained as subgraphs
  of $\Gamma$ since they have a vertex of degree~$3$.  Repeating this argument for all the $n$ vertices instead of
  vertex $1$ proves the claim for $n=3k$.
  
  Finally, let $n=3k+2$.  Again we use a similar argument as in the case $n=3k+1$ to get $3^k$ facets.  The
  corresponding graph $\Gamma$ has two edges not contained in any $3$-cycle.  Assume that one edge contains the
  vertex~$i$ and the other contains the vertex~$j$.  Consider $w\in\RR^n$ with $w_i=2$ and all other components equal to
  $1$.  We proceed as in the case $n=3k$, and again we apply Lemma~\ref{lem:B11}: As before the case
  \eqref{it:B11-cycle} is impossible because this would yield a spanning subgraph with at least $2+2(k-1)=2k<\lfloor
  2n/3\rfloor$ edges.  Hence we are in case \eqref{it:B11-non_cycle} to get a graph $\Gamma'$ with subgraphs of size
  $3+2(k-1)=2k-1=\lfloor 2n/3\rfloor$.  There are $3^{k-1}$ of that kind which are different from the spanning subgraphs
  of $\Gamma$ because $i$ has degree $2$.  A similar argument with $j$ instead of $i$ completes the proof of the
  theorem.
\end{proof}

We can also construct a metric for which this bound is tight. For an arbitrary graph $\Gamma$ on $n$ points we define a metric via
\[
d_\Gamma(i,j) \ = \
\begin{cases} 2;& \text{if } \{i,j\}\in E(\Gamma)\\
  1+\frac{1} {n^2+in+j} & \text{otherwise,}
\end{cases}
\]
for $i<j$.
Notice that our metric $\dmax^n$ corresponds to the graph on $n$ vertices without any edges.
We define $\dmin^n:=d_{\Gamma_{\min}^n}$ by letting
\begin{equation*}
  \{i,j\}\in E(\Gamma^n_{\min}) :\Leftrightarrow
  \begin{cases}
    \lfloor\frac{i-1}3\rfloor= \big\lfloor\frac{j-1}3\big\rfloor                       & \text{for $n\equiv 0,1\mod 3$}\\ 
    \lfloor{\frac{i-1}3}\rfloor= \big\lfloor{\frac{j-1}3}\big\rfloor \text{ and }i,j<n & \text{for $n\equiv 2\mod 3$} \quad .
  \end{cases}
\end{equation*}
So $\Gamma^n_{\min}$ consists of $\lfloor n/3\rfloor$ triangles and $n \mod 3$ isolated vertices. In fact, $d_{\min}^n$ is a slight
modification and generalization of the metric given by Develin in \cite[Proposition 3.3]{Develin} to proof the tightness of his lower
bound. Actually the proof that our bound is tight is obtained by analyzing the proof to  \cite[Proposition 3.3]{Develin} and
refining its techniques.

\begin{figure}[htbp]\centering
  \includegraphics[height=.35\textwidth]{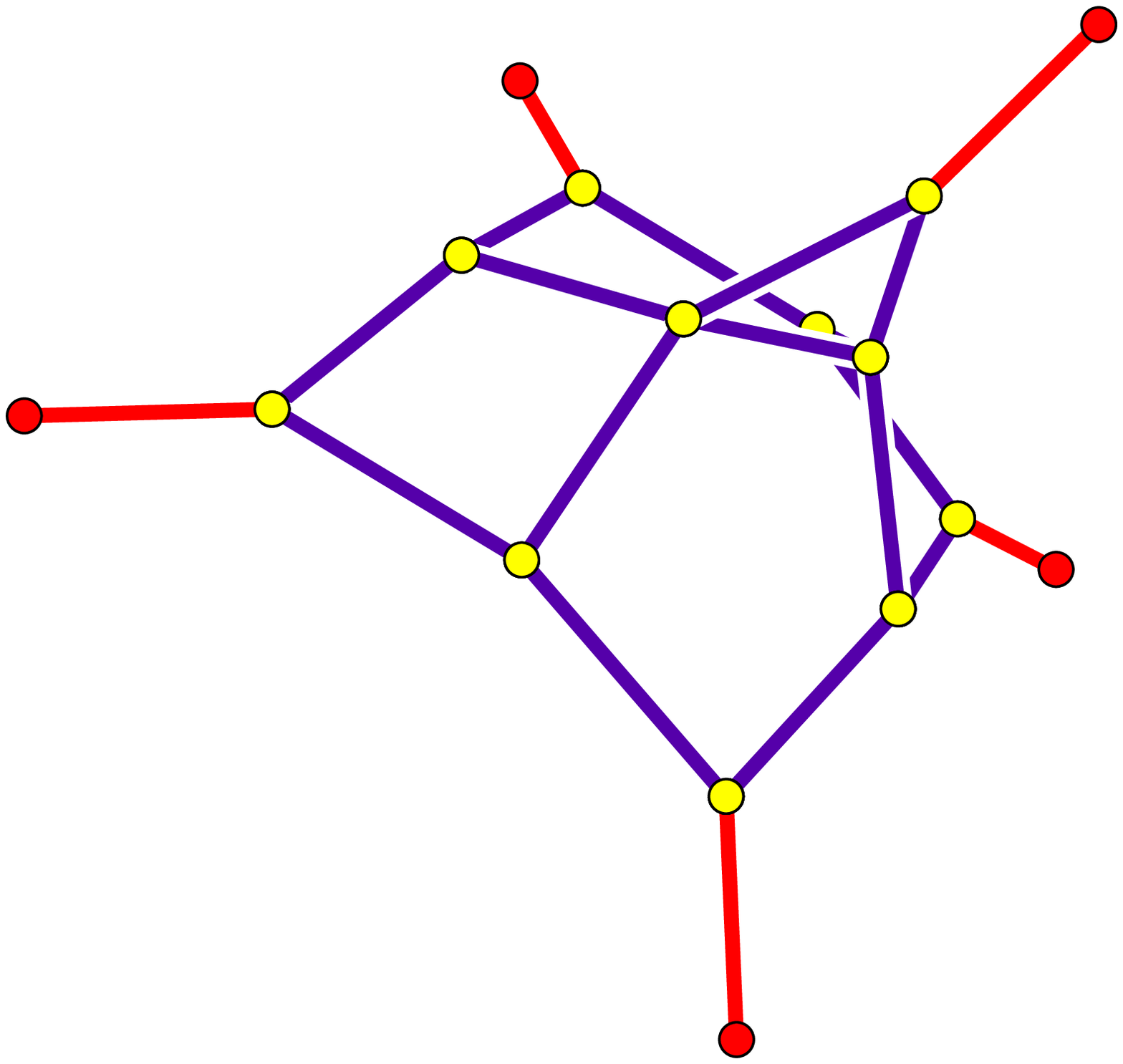}\quad
  \includegraphics[height=.35\textwidth]{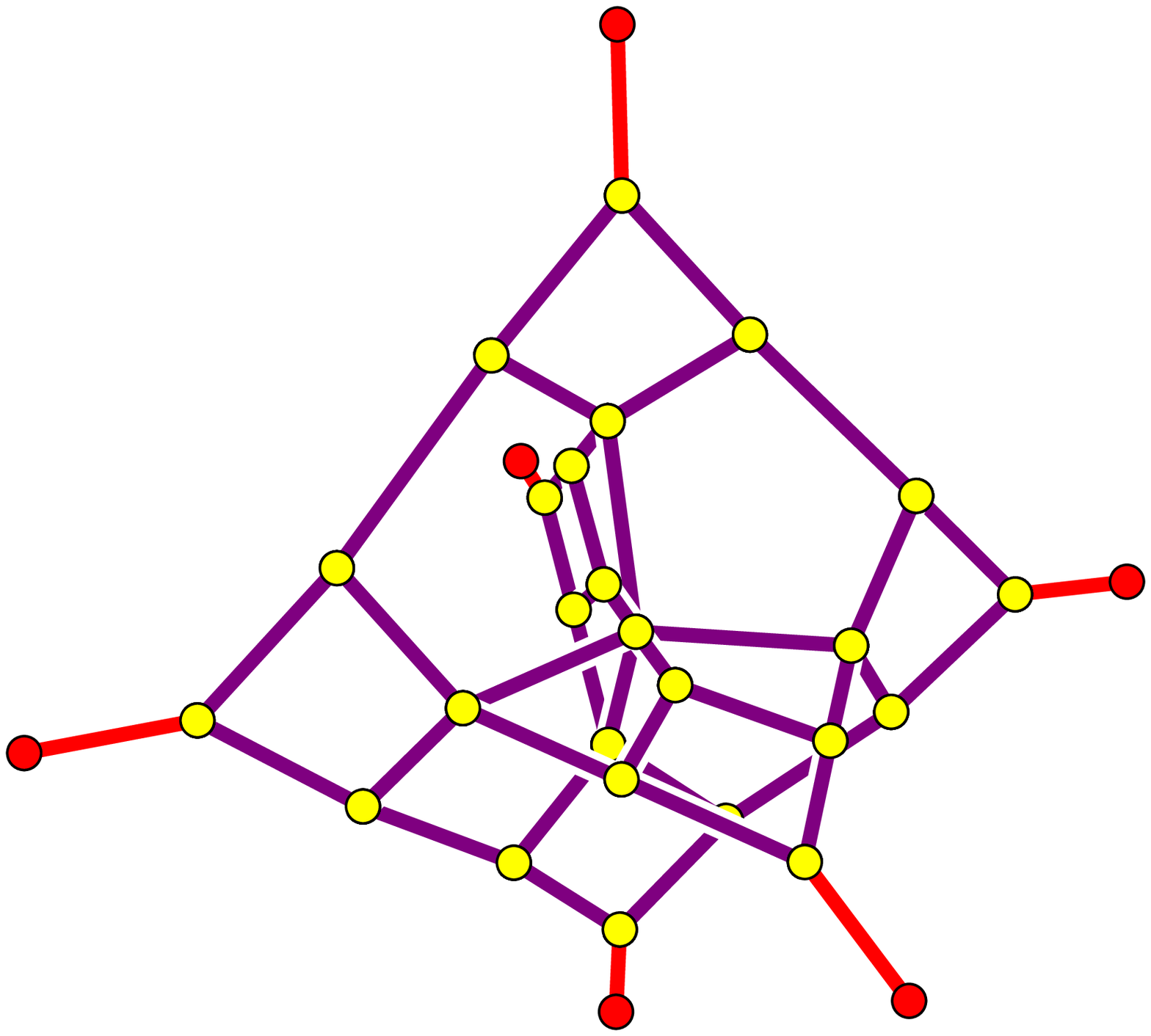}

  \caption{Visualization of (the graphs of) the tight spans $\Tmin^5$, with $f$-vector $(16,20,5)$, and $\Tmin^6$, with $f$-vector
    $(31,45,15)$.  Note that the image of $\Tmin^6$ shown is slightly misleading as the three collinear vertices in
    the center, in fact, form a triangle.  The tight span $\Tmin^6$, or rather the tight span of an equivalent ideal
    metric, occurs as $\#7$ in the list of Sturmfels and Yu~\cite{MR2097310}.}
\end{figure}

It is natural to ask if we can find a lower bound for all components of the $f$-vector from
Theorem~\ref{thm:lower-bound} in the same way as we derived Theorem \ref{thm:upper-bound} from Lemma
\ref{lem:n/2-bound}.  Unfortunately, this requires a much greater effort. The main problem is that there are
non-isomorphic subgraphs of $\Gamma_{\min}^n$ induced by submetrics of $\dmin^n$ of the same number of points; they even
give tight spans with different $f$-vectors.  Actually, such a proof would include the computation of the full
$f$-vector of all metrics $d_\Gamma$ with all components of $\Gamma$ of size at most~$3$.  Therefore, we suggest to
investigate the combinatorics and the $f$-vectors of the metrics $d_\Gamma$ for arbitrary graphs $\Gamma$.  This should
lead to a complete classification of all possible $f$-vectors of tight spans of generic metrics.
 
\bibliographystyle{amsplain}
\bibliography{main}

\end{document}